%% file: agt-5-37.tex
\documentclass{gtart_h}
\input agtout

\lognumber{37}
\volumenumber{5}
\volumeyear{2005}
\papernumber{37}
\pagenumbers{899}{910}
\received{10 February 2005} 
\accepted{1 June 2005}
\published{29 July 2005}

\usepackage{amssymb}

\usepackage[leqno]{amsmath}    

\numberwithin{equation}{section}        
\newtheorem{thm}[equation]{Theorem}     
\newtheorem{lem}[equation]{Lemma}       
\newtheorem{prop}[equation]{Proposition}   
\newtheorem{cor}[equation]{Corollary}      

\theoremstyle{definition}
\newtheorem*{D2}{D(2)--problem}
\def\vbar{\ \big|\ }
\def\mapright#1{\smash{\mathop{\longrightarrow}\limits^{#1}}}
\def\longlongrightarrow{\relbar\joinrel\!\relbar\joinrel\rightarrow}
\def\longmapright#1{\overset{#1}{\longlongrightarrow}}   
\def\ker{ker \hskip .25em}
\def\im{im \hskip .25em}
\def\mod{\hskip .4em mod \hskip .3em}
\def\G{{\mathbb Z}Q_{28}}
\def\Gn{{\mathbb Z}Q_{28}/\langle x^7+1\rangle }
\def\Gp{{\mathbb Z}Q_{28}/\langle \phi_{14}\rangle }
\def\Gi{{\mathbb Z}Q_{28}/\langle x+1\rangle }

\def\Q{{\mathbb Z}Q_{4n}}
\def\Qn{{\mathbb Z}Q_{4n}/\langle x^n+1\rangle }
\def\Qp{{\mathbb Z}Q_{4n}/\langle \psi_{2n}\rangle }
\def\Qi{{\mathbb Z}Q_{4n}/\langle x+1\rangle }

\begin{document}
\title[A stably free nonfree module]{A stably free nonfree module and 
its\\relevance for homotopy classification, case $Q_{28}$}
\asciititle{A stably free nonfree module and 
its\\relevance for homotopy classification, case Q_{28}}
\authors{F. Rudolf Beyl\\Nancy Waller }
\address{Department of Mathematics and Statistics, Portland 
State University\\Portland, OR 97207-0751, USA}

\email{beylf at pdx dot edu }

\begin{abstract} 
The paper constructs an ``exotic'' algebraic
\hbox{2--com}\-plex over the generalized qua\-ter\-nion
group  of order $28$, with the boundary maps given by
explicit matrices over the group ring.  This result depends on
showing that a certain ideal of the group ring is stably free
but not free. As it is not known whether the complex
constructed here is geometrically realizable, this example is
proposed as a suitable test object in the investigation of an
open problem of C.T.C.\ Wall, now referred to as the D(2)--problem.
\end{abstract}

\asciiabstract{%
The paper constructs an `exotic' algebraic 2-complex over the
generalized quaternion group of order 28, with the boundary maps given
by explicit matrices over the group ring.  This result depends on
showing that a certain ideal of the group ring is stably free but not
free. As it is not known whether the complex constructed here is
geometrically realizable, this example is proposed as a suitable test
object in the investigation of an open problem of C.T.C. Wall, now
referred to as the D(2)-problem.}



\asciikeywords{Algebraic 2-complex, Wall's D(2)-problem,
geometric realization of algebraic 2-complexes,
homotopy classification of 2-complexes, generalized quaternion groups,
partial projective resolution, stably free nonfree module}

\maketitle

The main topic of this paper is the construction of an ``exotic''
algebraic \hbox{2--com}\-plex over $Q_{28}$, the generalized
quaternion group of order $28$. The result provides significantly
more detail than mere existence proofs, as the boundary maps are
given by explicit matrices over the group ring (Section~\ref{Constr}).
This example should serve as a suitable test object in the investigation
of an open problem of C.T.C.\ Wall \cite[p.57]{Wa65}, now
\cite{JoB03} referred to as the D(2)--problem.

\begin{D2}\label{D2}
Suppose $X$ is a finite three-dimensional connected CW--com\-plex
(with universal cover $\widetilde X$\/) such that
$H_3(\widetilde X, \mathbb Z) = 0$ and $H^3(X, {\mathcal B})=0$
for all local coefficient systems ${\mathcal B}$ on $X$.
Is $X$ homotopy equivalent to a finite \hbox{2--com}\-plex?
\end{D2}

F.E.A.~Johnson \cite{JoS03}, \cite{JoB03} has shown that for finite
base groups this question has an affirmative answer if, and only if,
every algebraic \hbox{2--com}\-plex is geometrically realizable.

The work of K.W.~Gruenberg and P.A.~Linnell \cite[(2.6)]{GL88}
on minimal resolutions of lattices shows that stably equivalent
lattices over ${\mathbb Z}Q_{32}$, despite stringent minimality conditions,
need not have the same number of generators.  These lattices
can be viewed as the second homology modules of algebraic
\hbox{2--com}\-plexes. Johnson \cite{JoM03}, \cite{JoB03} shows
the existence of algebraic \hbox{2--com}\-plexes over $Q_{4n}$,
$ n \geq 6$, with ``exotic'' second homology modules and emphasizes
their relevance for the D(2)--problem.  It is not known whether these
examples are geometrically realizable. It is difficult to determine
whether these cases are actual counterexamples to C.T.C.\ Wall's
D(2)--problem since the maps and the free bases of the relevant modules
are not explicitly given.

Our joint work with M. Paul Latiolais \cite{BLW97} indicated to us
that an explicit potential counterexample could be constructed over
a binary polyhedral group $G$ which allows stably free nonfree
${\mathbb Z}G$--modules. Our example and those listed above depend on
Swan's work \cite{Sw83} on the existence of such projectives.
We give a specific example for $Q_{28}$ in Section~\ref{SFNF}.

The particular construction of an algebraic \hbox{2--com}\-plex
presented here grew out of a series of talks and workshops given
by us on geometrical realizability in Luttach (1997, 2001),
Chelyabinsk (1999), and Portland (2002). We would like to thank
Micheal Dyer, Cameron Gordon, Jens Harlander, Cynthia
Hog-Angeloni, Paul Latiolais, and Wolfgang Metzler for questions
and discussions that clarified the technical difficulties in
developing an explicit example.  Several of these contacts
were facilitated by NSF Grant INT--9603282 and a companion grant
by the DAAD for U.S.--Germany Cooperative Research.

\section{Preliminaries}\label{Prelims}

Given a group G, an {\em algebraic \hbox{2--com}\-plex\/} over
${\mathbb Z}G$ is a partial projective resolution of ${\mathbb Z}$,
$ C_2 \to C_1 \to C_0 \to {\mathbb Z} \to 0$, where the $C_i$ are
finitely generated, stably free ${\mathbb Z}G$--modules. We will usually
write this as an exact sequence of (left) ${\mathbb Z}G$--modules
\[
 A\co \qua 0 \to M \to C_2 \mapright{\partial_2} C_1
\mapright{\partial_1} C_0 \mapright{\partial_0} {\mathbb Z} \to 0\, ,
\]
where $M = \ker \partial_2 $ is the second homology module
and $M \to C_2 $ is inclusion. A familiar example of
an algebraic \hbox{2--com}\-plex is the cellular chain complex
$C(\widetilde K)$ of the universal cover $\widetilde K$ of a finite
two-dimensional CW--complex $K$ with fundamental group
$G$. Such a \hbox{2--com}\-plex is always homotopy equivalent
to some presentation complex for $G$.
An algebraic \hbox{2--com}\-plex which is chain homotopy
equivalent to the cellular chain complex of the universal cover of some
presentation complex for $G$ is called {\em geometrically realizable}.

One may also create new algebraic \hbox{2--com}\-plexes by
beginning with a known algebraic \hbox{2--com}\-plex $A$.  For the
reader's convenience, we summarize a construction due to Swan
\cite[\S6]{SwP60}. Suppose that there exists a ${\mathbb Z}G$--module
$M'$ such that $M' \oplus {\mathbb Z}G^k \cong M \oplus {\mathbb Z}G^k$
for some $k > 0$.

Form a new complex
\begin{equation}\label{1.1}
 0 \to M \oplus {\mathbb Z}G^k \longmapright{[i, id]}
 C_2 \oplus {\mathbb Z}G^k \to C_1 \to C_0 \to {\mathbb Z} \to 0 \, .
\end{equation}
Since ${\mathbb Z}G$ is weakly injective \cite{CE56}, \cite{SwP60},
there is another direct sum decomposition
$C_2 \oplus {\mathbb Z}G^k \cong Q_2 \oplus {\mathbb Z}G^k$ such that
\begin{equation}\label{1.2}
 0 \to M' \oplus {\mathbb Z}G^k \to Q_2 \oplus {\mathbb Z}G^k
  \to C_1 \to C_0 \to {\mathbb Z} \to 0
\end{equation}
is exact and the map restricted to ${\mathbb Z}G^k$ is the identity.
Note that $Q_2$ is stably free. The summand ${\mathbb Z}G^k$ may
now be cancelled to give
\[
 A'\co \qua 0 \to M' \to Q_2 \to C_1 \to C_0 \to {\mathbb Z} \to 0 \, .
\]
(To meet the requirement that $ M' \to Q_2 $ is inclusion, modify
$Q_2$ by renaming the elements in the image of $M'$.)
If the original complex $A$  is the chain complex of the
universal cover of a \hbox{2--com}\-plex $K$, then \eqref{1.1} corresponds
to $\overline K = K \vee\bigvee \limits _{i=1}^k S^2$. Since the
summand ${\mathbb Z}G^k$ to be cancelled in \eqref{1.2} represents elements of
$\pi_2(\overline K)$, attaching 3--balls to $\overline K$ gives a
\hbox{3--com}\-plex $K^3$ such that $C(\widetilde{ K^3})$ is chain
equivalent to $A'$. It is unknown whether $ A'$ is geometrically
realizable.

It is this question that remains unanswered in C.T.C.\ Wall's work
on the problem of determining whether a CW--complex is
homotopy equivalent to a finite one of lower dimension.

Consider the well-known presentation of the generalized
quaternion group of order $4n$:
$ Q_{4n} = \langle x, y \vbar x^n y^{-2}, y^2(xy)^{-2} \rangle$.
For the purposes of our construction, we
will use the following $Q^{**}$--equivalent presentation:
\[ Q_{4n} = \langle x,y \vbar x^ny^{-2}, yxyx^{-n+1}\rangle .\]
Let $K$ be the  associated presentation complex.
The relations of this presentation may be
represented by two disjoint simple closed curves on a genus 2
handlebody, and determine a Heegaard splitting of a 3--manifold
\cite[\S3.1]{HMS93} with fundamental group $Q_{4n}$. After an
appropriate thickening, $K$ may be
embedded in this manifold as a spine. Consequently $\pi_2(K) \cong
\Q/\langle N\rangle $, where $N \in  \Q$ is the sum of the group
elements. Schanuel's Lemma implies that any module  $M$, such that
$M\oplus \Q^k \cong \Q/\langle N\rangle \oplus \Q^k$ for some $k > 0$,
is a candidate for the second homology group of an algebraic
\hbox{2--com}\-plex  with minimal Euler characteristic over $\Q$.

Suppose $P$ is a $\Q$--projective with
\[
P \oplus \Q \cong \Q \oplus \Q \, ;
\]
that is, $P$ is {\em stably free}.  This isomorphism may be modified
by a basis change to give $(NP,0) \mapsto (N,0)$. Thus
\[
P/NP \oplus \Q \cong  \Q/\langle N\rangle  \oplus \Q \, .
\]
The construction earlier in this section yields an
algebraic \hbox{2--com}\-plex with second homology module
$P/NP$ and having the same Euler characteristic as that of a
3--manifold spine.

If $P/NP$ is not isomorphic to $\Q/\langle N\rangle $, this algebraic
\hbox{2--com}\-plex is
either chain homotopy equivalent to the cellular chain complex
of the universal cover of an {\em exotic\/} presentation complex
for $Q_{4n}$ ( i.e., not  homotopy equivalent to a
3--manifold spine), or is a counterexample to Wall's D(2)--problem.
This motivates the search for a suitable example of a stably free
nonfree module.

\section{A class of projective modules over $\Q$}\label{Class}

Given integers $a$ and $b$, let $P$ be the left ideal in $\Q$  generated
by $a+by$ and $x+1$. Since $(x+1)y = yx^{-1}(x+1)$
and $(a+by)x = x(a+by) + by(1-x^{-1})(x+1)$, $P$ is actually a
two-sided ideal.

\begin{prop}\label{P2.1}
Let $k=a^2+b^2$ for $n$ odd, and $k=a^2-b^2$ for $n$ even.
If $(k,2n)=1$ then $P$ is projective.
\end{prop}

\begin{proof}
Note that $k \in P$. Let $d=gcd(a,b)$ and $t = k/d $\ .
Since another generating set for $P$ is $\{t,x+1, a+by\}$,
$\Q/P \cong {\mathbb Z}/t{\mathbb Z} \oplus {\mathbb Z}/d{\mathbb Z}$
as abelian groups.  Since the quotient has exponent relatively prime
 to $4n$, $P$ is projective by \cite[Prop. 7.1, p.570]{SwI60}.
\end{proof}

Towards an explicit example of a stably free module that is not free,
from now on we assume $n$ odd and later will specialize to $n=7$.

Let $n$ be odd. In the spirit of Swan \cite[Ch.~10]{Sw83}, the ring
$\Q$ may be decomposed by a series of Milnor squares
\cite[\S2]{Mi71}, \cite[\S42]{CR87}:
\begin{equation}\tag\textup{I}
 \begin{array}{clc}
  \Q & \longrightarrow & \Q/\langle x^n+1\rangle \\
  \downarrow &    & \downarrow \\
  {\mathbb Z}D_{2n}  & \longlongrightarrow & {\mathbb F}_2 D_{2n}
 \end{array}
\end{equation}
\begin{equation}\tag\textup{II}
 \begin{array}{ccc}
  \Q/\langle x^n+1\rangle & \to & \Qp \\
  \downarrow &   & \downarrow \\
  \Qi  & \to & {\mathbb Z}_n[y]/\langle y^2+1\rangle
 \end{array}
\end{equation}
where $D_{2n}$ is the dihedral group of order $2n$ and
$\psi_{2n} = 1-x+x^2-x^3+\dots+x^{n-1}$.
Note that for an odd prime $p$, ${\mathbb Z}_p[y]/\langle y^2+1\rangle $
is isomorphic to the Galois field ${\mathbb F}_{p^2}$ of order $p^2$
when $p\equiv 3 \mod 4$, and isomorphic to
${\mathbb F}_p \times {\mathbb F}_p$ when $p \equiv 1 \mod 4$.
Every $\Q$--projective is isomorphic to some pullback
\[P(\, P_1,P_2, \alpha \,) = \{(a,b) \in P_1 \oplus P_2
\vbar \overline{a}= \alpha \overline{b} \in {\mathbb F}_2D_{2n} \},\]
where $P_1$ and $P_2$ are projectives over $\Qn$ and
${\mathbb Z}D_{2n}$, respectively, and $\alpha  \in Aut({\mathbb F}_2D_{2n})$,
cf.~\cite[(42.11)]{CR87}. Similarly, every $\Qn$--projective is isomorphic
to a pullback constructed from projectives over the quotient rings
in square~\textup{(II)}. Thus projective modules over $\Q$ are built (up to
isomorphism) by a sequence of pullbacks involving  related
projectives over quotient rings. For certain odd primes $p$, Swan
\cite[Ch.~10]{Sw83} has shown that suitable choices in this building
process lead to a ${\mathbb Z}Q_{4p}$--projective that is stably free
and not free.

\begin{prop}\label{P2.2}
Let $n$ be odd and let $(a^2 +b^2, 2n)=1$.
Let $P $ be the left ideal $\langle a+by, x+1\rangle $ in $\Q$.
Then $P/\langle x^n+1\rangle $  is isomorphic to
\begin{multline*}
\overline P = P(\, \Qp, \,\Qi, \,a+by \,) \\
 = \{(e,f) \in \Qp \oplus \Qi \vbar \overline{e}=(a+by)\overline{f}
\in {\mathbb Z}_n[y]/\langle y^2+1\rangle  \}
\end{multline*}
associated with square~\textup{(II)}, and $P$ is isomorphic to
$$P' = P(\,P/\langle x^n+1\rangle,\,{\mathbb Z}D_{2n},\,1\,) =
\{(c,d) \in P/\langle x^n+1\rangle  \oplus
{\mathbb Z}D_{2n} \vbar \overline{c}=\overline{d} \in {\mathbb F}_2D_{2n} \}$$
associated with square~\textup{(I)}.
\end{prop}

\begin{proof}
Consider square~\textup{(II)}.
Note that $a+by$ is a unit in ${\mathbb Z}_n[y]/\langle y^2+1\rangle $.
Define  $\eta_1\co P/\langle x^n+1\rangle  \to \Qp$ by the
projection $\Qn \to \Qp$ restricted to $P/\langle x^n+1\rangle $.
Any element $\rho=r(a+by) +s(x+1)$ of $P$ uniquely
determines  $r \mod \langle x+1\rangle $, because $(a+by)$ is not
a zero divisor in $\Qi$, the Gaussian integers. Thus $\eta_2(\rho)
= r + \langle x+1\rangle $ gives a well-defined homomorphism
$\eta_2\co P/\langle x^n+1\rangle  \to \Qi$. Since
$\overline{\eta_1(\rho)} =(a+by)\overline{\eta_2(\rho)}$ in
${\mathbb Z}_n[y]/\langle y^2+1\rangle $, this pair of maps uniquely
determines a homomorphism $\eta\co P/\langle x^n+1\rangle  \to
\overline P$ by the pullback property.

 A generating set for $\overline P$ in $\Qp \oplus \Qi $ is given by
\[ \{( a+by,1), (0, \psi_{2n}) , (x+1, 0), (0,n) , (n,0)\}. \]
Since $n\equiv \psi_{2n} \mod \langle x+1\rangle $, and $(0,
\psi_{2n})= \psi_{2n}(a+by, 1) $, this set may be reduced to the
two generators $(a+by,1), (x+1,0)$ which are the $\eta$--images of
$a+by,x+1$. Thus $\eta =(\eta_1,\eta_2)$ is an epimorphism. Since
${\mathbb Z}_n[y]/\langle y^2+1\rangle $ is finite, the rank of
$\overline P$ as a free abelian group is $2n$.  Since $P/\langle
x^n+1\rangle $ has the same rank, $\eta$ is an isomorphism.

Consider square~\textup{(I)}. Since one of $a$ or $b$ is even and the
other is odd, $P/\langle x^n+1\rangle $ maps onto ${\mathbb F}_2D_{2n}$.
Restricting the top homomorphism of square~\textup{(I)} to
$ P \to P/\langle x^n+1\rangle $
gives a commutative diagram like square~\textup{(I)} with $\Q$ replaced by
$P$. This defines a homomorphism $\tau\co P \to P'$, the pullback.
Now $P'$ is generated by
\[
\{(a+by, a+by), (x+1,x+1), (0,2)\}  \subseteq
P/\langle x^n+1\rangle \oplus {\mathbb Z}D_{2n}.
\]
The first two generators are obviously in $\im \tau$, and
$(0,2) = \tau(x^n+1)$. As before, the epimorphism $\tau$ is an
isomorphism, since $P$ and $P'$ have the same rank as abelian groups.
\end{proof}

\begin{cor}\label{C2.3}
If $\overline P$ is not free over $\Qn$, then $P$ is not free over $\Q$.
\end{cor}

\proof
As ${\langle x^n+1\rangle}P = \langle x^n+1\rangle $, this follows from
$$\overline P \cong P/{\langle x^n+1\rangle} \cong \Qn \otimes_{\Q} P\/ .
\eqno{\qed}$$

\begin{cor}\label{C2.4}
If $\overline P$ is not free over $\Qn$,
$P/NP$ is not isomorphic to $\Q/\langle N\rangle $.
\end{cor}

\begin{proof}
If $\overline P$ is not free over $\Qn$,
\[
\overline P \cong P/\langle x^n+1\rangle  \cong \Qn \otimes_{\Q} P/NP
\]
is not isomorphic to $\Qn \cong \Qn \otimes_{\Q} \Q/\langle N \rangle $.
\end{proof}

\section{A stably free module which is not free}\label{SFNF}

In this section we specialize to the case where $P$ is the left
ideal $\langle -3+4y,x+1\rangle $ of $\G$\/.

\begin{thm}\label{T3.1}
Let $P$ be the left ideal $\langle -3+4y,x+1\rangle$
of $\G$. Then $P$ is stably free and not free.
\end{thm}

The remainder of this section deals with the proof of this
theorem. Generalizing this result involves a number of special cases.
General criteria for a projective
$P = \langle a+by,x+1\rangle $ over $\Q$
to be stably free and not free are given in \cite{BW}.

\begin{thm}\label{T3.2}
Let $P$ be the left ideal $\langle -3+4y,x+1\rangle$ of $\G$.
Then $P$ is not free, and $P/NP$ is not isomorphic to $\G/\langle N\rangle $.
\end{thm}

\begin{proof}
Consider square (II) and associated pullbacks of the form
\begin{multline*}
 P(\,\Gp , \,\Gi , \,\alpha\,) \\
 = \{ (a,b) \in \Gp  \oplus \Gi \vbar
\overline{a}= \alpha \overline{b} \in \ {\mathbb F}_7[y]/\langle y^2+1\rangle  \},
\end{multline*}
where $\alpha$ varies over
$({\mathbb F}_7[y]/\langle y^2+1\rangle )^* \cong {\mathbb F}_{49}^*$.
In the proof of \cite[Lemma 10.13]{Sw83}, Swan shows that there
are four isomorphism classes which correspond to the cosets
$({\mathbb F}_7[y]/\langle y^2+1\rangle )^*/ \langle 3,y\rangle $,
with the free class represented by the trivial coset $[1]$.
These cosets form a cyclic group of order four generated
by the coset $[1+2y]$.
Since $[-3+4y] = [1+2y]^2$ is a nontrivial coset,
$\overline P$ is not free over $\Gn$. The assertion now
follows from Corollaries~\ref{C2.3} and  \ref{C2.4}.
\end{proof}

The work of Swan \cite[pp. 110--111]{Sw83} can be interpreted
to give that the set of stably free modules over $\G$
corresponds to the subgroup of order two in the above
construction.  While this approach leads to a proof that
$P$ is stably free, the following alternative proof  has the
advantage of giving an explicit basis for $\G \oplus P$ as a
free module. A computational method for constructing similar
isomorphisms over $\Q$ is given in \cite{BW}.

\begin{thm}\label{T3.3}
Let $P$ be the left ideal $\langle -3+4y,x+1\rangle $ of $\G$.
As a submodule of $\G \oplus \G$, $\G \oplus P$ is free with basis
$\{(\Phi_{11},\Phi_{12}),(\Phi_{21},\Phi_{22})\}$, where
\begin{align*}
\Phi_{11} &= x^7[1+(1-x^7)y][1+x^{-5}]-[7-7x^7-\Sigma^-][1+(x^{-3} +x^3)x^5y],\\
\Phi_{12} &= [1+(1-x^7)y][1-(x^{-3}+x^3)y]+[7-7x^7][1+x^5],\\
\Phi_{21} &= x^7[7\!-\!7x^7][1\!+\!x^{-5}]
- [19\!-\!20x^7][1\!-\!(1\!-\!x^7)y][1\!+\!(x^{-3}\!+\!x^3)x^5y] + 14\Sigma^-,\\
\Phi_{22} &= [7-7x^7-\Sigma^-][1-(x^{-3}+x^3)y] +[19-20x^7][1-(1-x^7)y][1+x^5],
\end{align*}
and $\Sigma^- = 1 -x+x^2-x^3+\dots +x^{12}-x^{13}$. Thus $P$ is stably free.
\end{thm}

\begin{proof}
Consider the map $\Phi\co \G \oplus \G \to \G \oplus \G$ given by \\
multiplication of row vectors on the right by the matrix
$\left[\begin{smallmatrix}
 \Phi_{11} & \Phi_{12}\\
 \Phi_{21} & \Phi_{22}
 \end{smallmatrix}\right]$.
In the rational group ring, this matrix factors as
\begin{multline*}
\begin{bmatrix}
 1+(1-x^7)y & 7-7x^7-\Sigma^- \\
 7-7x^7-\Sigma^- & (19-20x^7)[1-(1-x^7)y]
 \end{bmatrix} \times \\[7pt]
\begin{bmatrix}
 -x^7(1+x^{-5}) & 1-(x^{-3}+x^3)y \\
 1+(x^{-3}+x^3)x^5 y & 1+ x^5
\end{bmatrix}
\begin{bmatrix}
 \frac{14}{195} \Sigma^-\!{-}1 & 0 \\
 0 & 1
\end{bmatrix}.
\end{multline*}

Multiplying this product on the left first by
\begin{multline*}
 \begin{bmatrix}
  -1-(x^{-3}+x^3)y & -x^7(1+x^{-5}) \\
  -(1+x^5) & 1-(x^{-3}+x^3)x^5y
 \end{bmatrix} \times \\[7pt]
 \begin{bmatrix}
  (19-20x^7)(1-(1-x^7)y) & \Sigma^- -7 +7x^7 \\
  \Sigma ^- -7+7x^7 & 1+(1-x^7)y
 \end{bmatrix},
\end{multline*}
and then by
\[ (x+1)(x^2-x^7+x^{12}), \]
gives $ (x+1) \left[\begin{smallmatrix} 0 & 1 \\ 1 & 0 \end{smallmatrix}\right]$.
This shows that both
$ \begin{bmatrix} x+1 & 0 \end{bmatrix}$
and $\begin{bmatrix} 0 & x+1 \end{bmatrix} $
are in the image of $\Phi$.
Then  $\im \Phi \mod \langle x+1\rangle $ is generated by
\begin{gather*}
 \begin{bmatrix} 0 & (1+2y)(1+2y) \end{bmatrix}
  \equiv
 \begin{bmatrix}  0 & -3 +4y \end{bmatrix},\\
 \begin{bmatrix} -39(1-2y)(1+2y) +(14)(14) & 0 \end{bmatrix}
  \equiv
 \begin{bmatrix} 1 & 0 \end{bmatrix}.
\end{gather*}
Thus $\im \Phi = \G \oplus P \subseteq \G \oplus \G$. Since
$\G \oplus P$ and $\G \oplus \G$ have the same rank as free
abelian groups,  $\Phi $ restricts to an isomorphism.
\end{proof}

\section{The construction}\label{Constr}

Let $n$ be odd and let $(a^2 +b^2, 2n) =1$. Let $P $ be the left
ideal $\langle a+by, x+1\rangle $ in $\Q$. The following
construction gives a partial projective resolution which is an
algebraic 2--complex whenever $P$ is stably free.

Let $K$ be the presentation complex associated with
\[ Q_{4n} = \langle x,y \vbar x^ny^{-2}, yxyx^{-n+1}\rangle\, .\]
Denote the relators by $R_1 = x^ny^{-2}$ and $R_2 = yxyx^{-n+1}$.
Let $C(\widetilde K)$ be the cellular chain complex  of the
universal cover of $K$:

\begin{equation}\label{4.1}
 0 \to H_2(\widetilde K) \to \Q \oplus \Q \mapright{\partial _2}
 \Q \oplus \Q \mapright{\partial _1} \Q \to {\mathbb Z} \to 0 \, .
\end{equation}

\begin{lem}\label{L4.2}
$ H_2(\widetilde K)$ is generated by
$$ \sigma = (1-y) \widetilde R_1 + (1-yx)\widetilde R_2 \, ,$$
and isomorphic to $ \pi_2(K) \cong \Q/\langle N\rangle $
where $ 1 + \langle N \rangle \mapsto \sigma $.
This isomorphism induces  $P/NP \cong P\sigma$.
\end{lem}

The authors originally found the formula for $\sigma$ by inspecting
the Heegaard splitting determined by these relators.
Here is an alternative algebraic proof.

\begin{proof}
By \cite[Prop. 1.2]{BLW97},
$H_2(\widetilde K) \cong \pi_2(K) \cong \Q/\langle N\rangle $.
By inspection,\\
$\sigma \in H_2(\widetilde K) \, $.
If $\sigma$ does not generate $H_2(\widetilde K)$,
then there is a nonzero solution to
\begin{align*}
A[1+x+x^2+\dots+x^{n-1}] +B[y-1-x-x^2-\dots-x^{n-2}]&= 0 \\
                  A[-y-1] + B[1+yx] &=0
\end{align*}
with $B= B_1(x) + B_2(x)y$ and $A$ of the form $A=A_1(x)$.
Comparing coefficients with respect to $y$ gives
\begin{align*}
A_1(x)[1\!+\!x\!+\!x^2\!+\dots+\!x^{n-1}] +
B_1(x)[-\!1\!-\!x\!-\!x^2-\dots-\!x^{n-2}] + B_2(x)x^n &= 0 \\
  B_1(x) - B_2(x)[1+x^{-1}+x^{-2}\!+\dots + x^{-n+2}]&=0\\
                     -A_1(x) + B_1(x) + B_2(x)x^{n+1} &=0 \\
                       -A_1(x) + B_1(x)x^{-1} +B_2(x) &=0 .
\end{align*}
Every solution of this system is a multiple of
$$ (1-x^n)\widetilde R_1 + (1 - x^{n+1} + [1-x^{-1}]y)\widetilde R_2 =
(1+y)\sigma\, .$$
Since $\langle a+by, x+1\rangle $ always contains an element with
augmentation $1$, $NP = \langle N\rangle $.
Thus the kernel of the composite map
$$ P \hookrightarrow \Q \twoheadrightarrow \Q/\langle N\rangle
\cong \langle \sigma\rangle $$
is $\langle N\rangle $ and the image is $P\sigma$.
Therefore $P/NP \cong P\sigma$.
\end{proof}

Note that  $P\sigma$ is generated by
\begin{align*}
 (x+1)\sigma &= (x+1)[(1-y)\widetilde R_1 + (1-yx)\widetilde R_2] \\
 \hbox{\rm and}\quad\quad\quad
  (a+by)\sigma &= (a+by)[(1-y)\widetilde R_1 + (1-yx)\widetilde R_2]\, .
\end{align*}
Consider the subcomplex of~\eqref{4.1} given by
\begin{equation}\label{4.3}
 0 \to \ker \partial'_2 \to \Q \oplus P \mapright{\partial'_2}
\Q \oplus \Q \mapright{\partial_1} \Q \to {\mathbb Z} \to 0 \,,
\end{equation}
where $\partial'_2$ is the restriction.

\begin{prop}\label{P4.4}
The kernel of $\partial'_2 $ equals $P\sigma$ and the sequence~\eqref{4.3} is
exact.
\end{prop}

Proposition~\ref{P4.4} and Lemma~\ref{L4.2} combine to imply that the
second homology module $\ker \partial'_2 $ of~\eqref{4.3}
is isomorphic to $P/NP$.

\begin{proof}
 (a)\qua To verify that $\ker \partial'_2 = P\sigma $:
Since $P$ is a two-sided ideal of $\Q$,
$P\sigma \subseteq \Q \oplus P \subseteq  C_2(\widetilde K)$.
Since $\ker \partial_2$  is generated by
\[
\sigma = (1-y)\widetilde R_1 + (1-yx)\widetilde R_2\, ,
\]
an element $C\sigma \in \ker \partial _2$ lies in $\Q \oplus P$
if, and only if, there exist $ A,B \in \Q$ with
$A(x+1) + B(a+by) = C(1-yx)$.  Since $n$ is odd, reducing
this problem$\mod \langle x+1\rangle $ gives the condition:
$$ (B_1 +B_2y)(a+by)  \equiv (C_1 + C_2y)(1+y) \text{ where }
  B_1,B_2, C_1,C_2 \in {\mathbb Z}\, .$$
Thus
\[
 2 \begin{bmatrix} C_1 & C_2 \end{bmatrix}
  =
 \begin{bmatrix} B_1 & B_2 \end{bmatrix}
 \begin{bmatrix}
  a+b & b-a \\
  a-b & a+b
  \end{bmatrix}.
\]
This equation$\mod 2$ implies $ B_2 - B_1 $ is even.
For an integral solution,
$ \begin{bmatrix} C_1 & C_2 \end{bmatrix} $
must be an integral linear combination of
$ \begin{bmatrix} a & b \end{bmatrix} $
 and
$ \begin{bmatrix} a-b & a+b \end{bmatrix} $.
Since
$(a-b)+(a+b)y \equiv (y+1)(a+by) \mod \langle x+1\rangle $,
$C$ is a $\Q$--linear combination of $a+by $ and $x+1$. Thus
$C\sigma \in \ker\partial'_2 $ implies $C \in P$ and
$C\sigma \in  P\sigma$.

(b)\qua To verify that $\im \partial'_2 = \im \partial_2 = \ker \partial_1$:
Observe that $C_2(\widetilde K) = \langle \sigma \rangle + (\Q \oplus P)$.
Since $\widetilde R_1 \in \Q \oplus P $, it is only necessary to show that
$ \widetilde R_2 \in \langle \sigma \rangle + (\Q \oplus P)$.
The latter submodule  contains the multiples of
$ \widetilde R_2 $ with the coefficients $a+by$, $x+1$, $1-yx$,
thus also with
\begin{align*}
 y+1 &= y(x+1) + (1-yx) ,\\
a-b &= (a+by)-b(y+1) ,\\
2 &= \psi_{2n}(x)\cdot (x+1) + (1-y)(y+1) ,
\end{align*}
and eventually $1 \cdot \widetilde R_2$ since $a-b$ is odd.
Since $\partial _2 \sigma = 0$, the result follows.
\end{proof}

Now specialize to $n = 7$ and
$P= \langle -3+4y, x+1\rangle \subseteq \G$\,.
By Theorems~\ref{T3.3} and  \ref{T3.2}, $P$ is stably free,
but $P/NP$ is not isomorphic to $\G/\langle N\rangle$\,.
In this case our complex is an {\em exotic\/} algebraic
\hbox{2--com}\-plex over $\G$\,. That is, either this complex is
geometrically realizable with fundamental group $\G$ and is
not homotopy equivalent to a 3--manifold spine, or is a
counterexample to Wall's D(2)--problem.
Viewing $\G \oplus P $ as a free module, with the basis specified
in Theorem~\ref{T3.3}, gives a computable example:

\begin{thm}\label{T4.5}
For $G = Q_{28}$ the following free algebraic \hbox{2--com}\-plex is
exotic:
$$ {\mathbb Z}G^2 \mapright{\partial_2} {\mathbb Z}G^2 \mapright{\partial_1}
{\mathbb Z}G \mapright{\varepsilon} {\mathbb Z} \to 0 \, , $$
where $\varepsilon$ is augmentation and the boundary maps are given
by multiplying \\ row vectors on the right by matrices.
The matrix for $\partial_1$ is
$ \begin{bmatrix} x-1 \\ y-1 \end{bmatrix} $
and the matrix for $\partial_2$ is
\[
 \begin{bmatrix}
  \Phi_{11}\sum \limits _{k=0}^6 x^k + \Phi_{12}(y-\sum \limits_{k=0}^5 x^k)
   &\quad -\Phi_{11}(1+y) + \Phi_{12}(1+yx) \\
  \Phi_{21}\sum \limits _{k=0}^6 x^k + \Phi_{22}(y-\sum \limits_{k=0}^5 x^k)
   &\quad -\Phi_{21}(1+y) + \Phi_{22}(1+yx)
 \end{bmatrix}
\]
where the $\Phi_{ij}$ are defined in Theorem~\ref{T3.3}.
\qed \end{thm}

\Addressesr
\end{document}

%% file: agtout.tex
\def\ifplaintex{\expandafter\ifx\csname documentclass\endcsname\relax}

\def\gtp{{\mathsurround=0pt\it $\cal G\mskip-2mu$eometry \&\ 
$\cal T\!\!$opology $\cal P\!$ublications}}  

\def\Addressesr{\bigskip
{\small \parskip 0pt \leftskip 0pt \rightskip 0pt plus 1fil \def\\{\par}
\sl\theaddress\par
\medskip
\rm Email:\stdspace\tt\theemail\hfill\rm Received:\qua\receiveddate \par}}

\def\recd{{\small Received:\qua\receiveddate\ifx\reviseddate\relax
\else\qquad Revised:\qua\reviseddate\fi\par}} 

\def\lognumber#1{\def\thelognumber{#1}}
\def\volumenumber#1{\def\thevolumenumber{#1}}
\def\volumeyear#1{\def\thevolumeyear{#1}}
\def\papernumber#1{\def\thepapernumber{#1}}
\def\pagenumbers#1#2{\def\startpage{#1}\def\finishpage{#2}}
\def\published#1{\def\publishdate{#1}}

\def\received#1{\def\receiveddate{#1}}

\def\accepted#1{\def\accepteddate{#1}}
\def\asciititle#1{\def\theasciititle{#1}}

\long\def\asciiabstract#1{\long\def\theasciiabstract{#1}}
\def\asciikeywords#1{\def\theasciikeywords{#1}}

\let\\\par\let\thelognumber\relax\let\thevolumenumber\relax
\let\thepapernumber\relax\let\thevolumeyear\relax\let\startpage\relax
\let\finishpage\relax\let\publishdate\relax\let\receiveddate\relax
\let\reviseddate\relax\let\accepteddate\relax\let\theasciititle\relax
\let\theasciiauthors\relax
\let\theasciiabstract\relax\let\theasciikeywords\relax

\let\theasciiemail\relax

\ifplaintex
\font\logobig=cmssbx10 scaled 3836
\font\logomed=cmssbx10 scaled 2557
\else
\font\logobig=cmssbx10 scaled 4200
\font\logomed=cmssbx10 scaled 2800
\fi

\long\def\makeagttitle{   
\count0=\startpage
\agt\hfill      
\hbox to 45truept{\vbox to 0pt{\vglue -13truept{\logomed A\kern -.37em{\logobig 
T}\kern -.38em G}\vss}\hss}
\break
{\small Volume \thevolumenumber\ (\thevolumeyear)
\startpage--\finishpage\nl
Published: \publishdate}

\vglue .25truein

{\parskip=0pt\leftskip 0pt plus
1fil\def\\{\par\smallskip}{\Large\bf\thetitle}\par\medskip} \vglue
0.05truein

{\parskip=0pt\leftskip 0pt plus 1fil\def\\{\par}{\sc\theauthors}
\par\medskip}%
 
\vglue 0.03truein 

{\small\leftskip 25truept\rightskip 25truept{\bf Abstract}\stdspace\theabstract

{\bf AMS Classification}\stdspace\theprimaryclass
\ifx\thesecondaryclass\relax\else; \thesecondaryclass\fi\par
{\bf Keywords}\stdspace \thekeywords\par}\vglue 7truept

}   

\ifplaintex
\font\phead=cmsl9 scaled 950
\font\pnum=cmbx10 scaled 913
\font\pfoot=cmsl9 scaled 950
\headline{\vbox to 0pt{\vskip -4.5mm\line{\small\phead\ifnum
\count0=\startpage ISSN 1472-2739 (on-line) 1472-2747 (printed)
\hfill {\pnum\folio}\else\ifodd\count0\def\\{ }%
\ifx\theshorttitle\relax\thetitle\else\theshorttitle\fi\hfill{\pnum\folio}
\else\def\\{ and }{\pnum\folio}\hfill\ifx\theshortauthors\relax\theauthors
\else\theshortauthors\fi\fi\fi}\vss}}
\footline{\vbox to 0pt{\vglue 0mm\line{\small\pfoot\ifnum\count0=\startpage
\copyright\ \gtp\hfill\else
\agt, Volume \thevolumenumber\ (\thevolumeyear)\hfill\fi}\vss}}
\else
\font=cmsl9 scaled 1050
\font\lnum=cmbx10 
\font=cmsl9 scaled 1050
\makeatletter
\def\@oddhead{{\small\lhead\ifnum\count0=\startpage ISSN 1472-2739 
(on-line) 1472-2747 (printed)\hfill {\lnum\number\count0}\else\ifodd\count0
\def\\{ }\ifx\theshorttitle\relax \thetitle \else\theshorttitle\fi\hfill
{\lnum\number\count0}\else\def\\{ and }{\lnum\number\count0}
\hfill\ifx\theshortauthors\relax 
\theauthors\else\theshortauthors\fi\fi\fi}}\def\@evenhead{\@oddhead}
\def\@oddfoot{\small\lfoot\ifnum\count0=\startpage\copyright\ \gtp\hfill\else
\agt, Volume \thevolumenumber\ (\thevolumeyear)\hfill\fi}
\def\@evenfoot{\@oddfoot}
\makeatother
\fi
\let\maketitlepage\makeagttitle
\let\makeshorttitle\maketitlepage
\let\maketitle\maketitlepage

\newwrite\gtoutfile
\long\gdef\makeheadfile{  
{\def\\{, }\def\s{ }
\immediate\openout\gtoutfile head.xxx
\immediate\write\gtoutfile{Proxy-for: \ifx\theasciiauthors\relax
\theauthors\else\theasciiauthors\fi\s<\ifx\theasciiemail\relax\theemail\else\theasciiemail\fi>}
\immediate\write\gtoutfile{\noexpand\\}
\immediate\write\gtoutfile{Authors: \ifx\theasciiauthors\relax
\theauthors\else\theasciiauthors\fi}
{\def\\{ }\immediate\write\gtoutfile{Title: \ifx\theasciititle\relax
\thetitle\else\theasciititle\fi}}
\immediate\write\gtoutfile{Subj-class: GT or SG, GR etc}
\immediate\write\gtoutfile{MSC-class: \theprimaryclass\ifx\thesecondaryclass\relax\else, \thesecondaryclass\fi}
\immediate\write\gtoutfile{Journal-ref: Algebr. Geom. Topol. \thevolumenumber\s
(\thevolumeyear) \startpage-\finishpage}
\immediate\write\gtoutfile{Comments: Published by Algebraic and
Geometric Topology at}
\immediate\write\gtoutfile{\s\s\s  http://www.maths.warwick.ac.uk/agt/AGTVol\thevolumenumber/agt-\thevolumenumber-\thepapernumber.abs.html}
\immediate\write\gtoutfile{\noexpand\\}
\immediate\write\gtoutfile{}
\ifx\theasciiabstract\relax
\immediate\write\gtoutfile{\theabstract}\else
\immediate\write\gtoutfile{\theasciiabstract}\fi
\immediate\write\gtoutfile{}
\immediate\write\gtoutfile{\noexpand\\}
\immediate\write\gtoutfile{}
\immediate\closeout\gtoutfile}}  

\def\maketitlepage{\makeagttitle\makeheadfile}
\let\makeshorttitle\maketitlepage
\let\maketitle\maketitlepage